%% file: main.tex
\documentclass[a4paper, 12pt]{article}
\usepackage[left=3cm,right=2cm,top=2.5cm,bottom=2.5cm]{geometry}

\usepackage[utf8]{inputenc}
\usepackage{lmodern,stackrel,verbatim,textcomp}
\usepackage[T1]{fontenc}
\usepackage[pdftex]{color,graphicx}
\usepackage{upgreek,graphicx, epstopdf,dsfont,soul}
\usepackage{mathtools,listings,marginnote}
\usepackage{amsmath,amsfonts,amssymb,amsthm,mathrsfs,stackrel}
\usepackage{xargs} 
\usepackage{tikz}
\usetikzlibrary{shapes,arrows,calc,positioning}
\usepackage{pgfplots}
\usepackage{xcolor}
\usepackage{caption}
\usepackage{subcaption}

\usepackage[normalem]{ulem}

\pgfplotsset{compat=1.17}

\usepackage[nottoc]{tocbibind}
\usepackage{hyperref}
\hypersetup{
    colorlinks=true,
    linkcolor=blue,
    filecolor=blue,
    urlcolor=blue,
    }

\RequirePackage[textwidth=16mm, textsize=tiny]{todonotes}

\renewcommand{\leq}{\leqslant}
\renewcommand{\geq}{\geqslant}

\newcommand{\ind}{\mathbf{1}}
\renewcommand{\P}[1]{\mathbb{P}\left(#1\right)}
\newcommand{\e}{\varepsilon}

\newcommand{\N}{{\mathbb{N}}}
\newcommand{\R}{\mathbb{R}}

\newcommand{\Ex}[1]{\mathbb{E}\left(#1\right)}

\newcommand{\Fr}{F_{\omega}}
\theoremstyle{plain}
\newtheorem{Thm}{Theorem}

\newtheorem{Prop}{Proposition}

\newtheorem{Lem}{Lemma}

\theoremstyle{definition}

\newtheorem{Rmk}{Remark}

\begin{document}

\title{Random walk in a rotational environment}
\author{Alberto M. Campos\thanks{Phd Student from IMPA- Instituto de Matematica Pura e Aplicada- RJ} \and Tarcísio P.R. Campos\thanks{Research collaborator at CDTN- Centro de Desenvolvimento de Tecnologia Nuclear, research collaborator at CNEN - Comissão de Energia Nucear, and voluntary research collaborator at UFMG- Univercidade Federal de Minas Gerais.}}
\date{\today}
\maketitle

\begin{abstract}
We define a random walk of a particle in $\R^3$ where the space is rotating. The particle is not glued to the space and will collide with it at random times, resulting in changes in its velocity and direction.  After many collisions, the random walk starts to have some asymptotic behaviors
inherited from the movement of space. The paper will find the limit movement of the particle, and explain how the randomness of the random walk gives rise to the particle asymptotic deterministic movement.
\end{abstract}

{\footnotesize Keywords: Langevin's equation ;random walk; Brownian Motion}
 
\input{1}

\input{2}

\input{3}

\bibliographystyle{plain}
\bibliography{ref}
\end{document}

%% file: 1.tex
\section{Introduction}\label{sec:Intro}\noindent

The first mention of the Brownian motion in the literature is related to the movement of pollen particles, where in still water the pollen perform a movement that approximates the Brownian motion; the original paper was not published, but one can follow the same experiment in \cite{BrownianExperiment}. Curiously, if the water moves,  the trace of the pollen particle is almost deterministic following the flow of the water. This raises an intriguing open question: How does the movement of the medium affect the random motion and is it possible to explain the deterministic component of its movement induced by the environment using its randomness?

The paper focuses on understanding how the medium affects the random motion in a specific case where the particle is moving in a rotational environment. In the problem formally stated in Section \ref{sec:Notation}, the particle is floating in the medium, moving along a straight line without interference from the space in which it is traveling until a certain point in time. At that moment, the particle with certain relative velocity will collide with the space, altering its direction and velocity.

Without solving mathematically due to the high number of collisions, one can stipulate a solution by a physical approximation. The particle as the process evolves will rotate following the velocity of the space, be dragged away by a centripetal force, and also have a random perturbation due to the random nature of the problem, see \cite{Rotationparticleslab} for some experiments and simulations of such movement. 

\begin{figure}[ht!]
    \centering
    \includegraphics[scale=0.5]{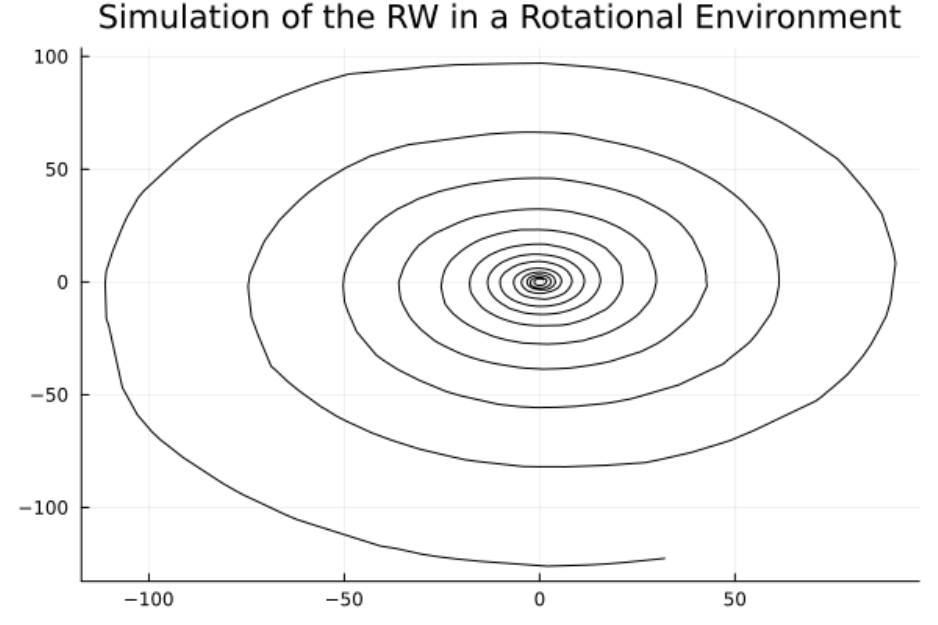}
    \caption{Simulation of the model of a random walk over a rotational environment introduced in Section \ref{sec:Notation}.}
    \label{fig:A}
\end{figure}

In this article, we will analyze the impact of multiple collisions over the movement, showing that the physics stipulation needs further analysis. Spoiling the result, the particle will follow the rotation movement of the space, but the centripetal force will not appear with the same intensity as expected. The particle is walking in a straight line except in collisions, and thus this force that drags the particle always from the origin is not a classic centripetal force, existing only on the average results of these random collisions.  This paper is focused on explaining this limit behavior and shows in a combination of random and geometric arguments why the particle behaves as it does; see Theorems \ref{Thm:1} and \ref{Thm:2}, and compare the process in Figure \ref{fig:A}. 

The model proposed in this paper shares similarities with others found in the literature that investigate random walks subjected to external forces; see \cite{Patlak,RandonWalkExternalForce1,RandonWalkExternalForce2,RandonWalkExternalForce3}. It is not surprising that the model exhibits these similarities, as it seeks to approximate a physical problem while preserving the natural laws of physics. Although similar in some topics, the paper is different from the previous literature in some points. The first point of difference is the physical formulation done in those articles, while here a purely mathematical approach will be exposed. The second point of difference is related to the fixed choice of motion, where by fixing the type of movement the result becomes more specific with simpler formulas and deeper results. The best example is the work \cite{Patlak} of Patlak that considers general types of movement, but does not find limits in distributions for the random walk. 

\begin{figure}[ht!]
    \centering
    \includegraphics[scale = 0.5]{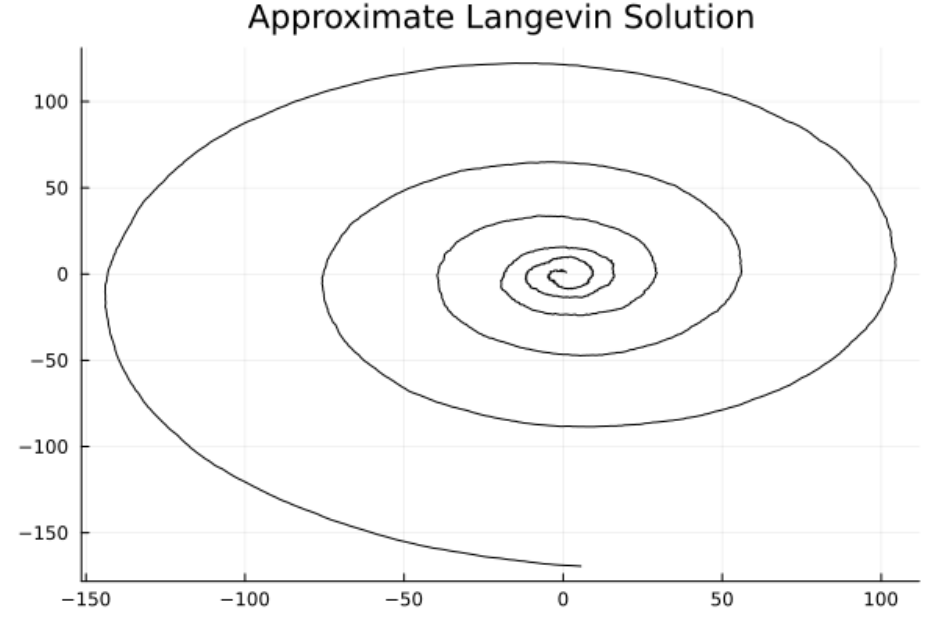}
    \caption{Simulation of the solution of the Lagevin equation \eqref{eq:Intro-Lagevin-Possible-Approximation}. }
    \label{fig:B}
\end{figure}

Perhaps the greatest instrument for understanding the behavior of a particle subject to random and deterministic forces is the Langevin equations; see \cite{Langevin-Intr-Apl-Book,Langevin-Book}. Although it is a great approximation for the expected movement of the problem, the method does not explain how these random sums of interactions converge to this deterministic behavior.  To illustrate those ideas, consider the following Lagevin's equation that follows the physical stipulation of the random walk in a rotational environment:
\begin{align}\label{eq:Intro-Lagevin-Possible-Approximation}
    X'(t)= \begin{bmatrix}
0 & -1 \\
1 & 0 
\end{bmatrix}  X(t) + \begin{bmatrix}
0.1 & 0 \\
0 & 0.1 
\end{bmatrix}  X(t) +dB_t.
\end{align} where $X(t)$ is the position of the particle, and $dB_t$ is a white noise. The first matrix in the equation is responsible to rotate the particle in the velocity of the space, while the second matrix is responsible to give it the associated centripetal force.  Notice that if the particle rotates in the velocity of the space a centripetal force exists, thus equation \eqref{eq:Intro-Lagevin-Possible-Approximation} is precisely what happens if the movement of the particle is continuously force to rotate by the space. Observe a simulation of the solution of the equation \eqref{eq:Intro-Lagevin-Possible-Approximation} in Figure \ref{fig:B}.

{}

Comparing the solutions in Figures \ref{fig:A} and \ref{fig:B} side by side, the force that drifts the particles presents different growth rates; this and other differences are treated in detail in Section \ref{sec:Notation}, by Remark \ref{Rmk:Phiscs}. Essentially the physics stipulation does not capture the real behavior of the mathematical model presented here. 

This paper was divided as follows. In Section \ref{sec:Intro} a brief introduction and motivation of the problem are given. The Section \ref{sec:Notation} exposed the definition of the model and the problems proved in this paper. Finally, in Section \ref{sec:rot} we proved the theorems for a random walk in the rotational environment. 

\textbf{Acknowledgments} Both authors thank Renato Soares dos Santos for his valuable comments and the CNPq Conselho Nacional de Desenvolvimento Científico Tecnolôgico. Research  A.M.C. is supported in part by CNPq, grant
141068/2020-5

%% file: 2.tex
 \section{Notations and main theorems}\label{sec:Notation}
\noindent 

Consider $\mathbb{R}^3$ as the three-dimensional space. The variables $x=(x_1,x_2,x_3)$ and $y=(y_1,y_2,y_3)$ belonging to $\mathbb{R}^3$ are used to denote points in the space and eventually will be interpreted as positions. Additionally, the variables $v=(v_1,v_2,v_3)$, $u=(u_1,u_2,u_3)$ and $e=(e_1,e_2,e_3)$ represent vectors in $\mathbb{R}^3$, and for us always will be interpreted as velocity.  

The variable $t$ is used exclusively to represent time and the variable $m$ will be used exclusively as the mass of a particle, where $t,m\in \mathbb{R}_+=\{z\in \mathbb{R} :z\geq 0\}$.

Define the $L^2$ norm of $\R^3$ as $\|\cdot\|_{}$, where for every $x\in \mathbb{R}^3$, $\|x\|=\sqrt{x_1^2+x_2^2+x_3^2}$. Also, define the open ball with radius $r>0$ and center $x\in \mathbb{R}^3$ as the set $B(x,r)=\{y\in \R^3: \|x-y\|<r\}$. In this paper, the partial distance taken at the first two coordinates in relation to third variable axis is relevant, so define:
\begin{align*}
    d_r(x)=\sqrt{x_1^2+x_2^2}, 
\end{align*} to be the radial distance.

In the paper, whenever some relevant quantity is random, we denote it by a capital letter. So, $X$ is used as a random position in space, $V$ is used as a random velocity, $E$ is used as an error vector of velocity, and $T$ is used as a random time.   

Define a particle at time $t>0$ as a triplet $(X_t,V_t,t)$, where $X_t\in \R^3$ is the position of the particle, $V_t\in \R^3$ is the velocity and $t$ is the time, also set the initial values $X_0=x_0$, $V_0=v_0$. The particle is not bound or glued to the space and will interact with it in the form of collisions at random times $(T_k)_{k\in \N}$. The choice of these times depends on some factors that will be explained later. As an abuse of notation, with the times $(T_k)_{k\in \N}$ fixed, denote for every $k\in \N$ the values $X_k=X_{T_k}$ and $V_{k}=V_{T_k}$ as the position and the velocity in the $k$-th interaction of the particle. In particular, for $t\in[T_k,T_{k+1}]$, the position of the particle is $X_t=X_{k}+(t-T_k)V_{{k}}$. 

The space in which the particle belongs also performs some independent movement. Since the particle floats over the space, the linear movement of the particle over $\mathbb{R}^3$ will be compounded by the movement of the space creating a new path (not necessarily a line),  denoted by $\gamma \subset \mathbb{R}^3$. Let $(\gamma_k)_k$ be the paths that the particle goes through between times of interactions $(T_k)_k$, and set $|\gamma|$ to be the length of any fixed path $\gamma$.  

In this paper, we restrict the movement of the environment to be a rotation. Since, for any rotation of $\R^3$ one axis is made fixed. Without loss of generality, consider, for any point $x\in \mathbb{R}^3$, $\omega>0$, and time $t>0$, the movement position: 
\begin{align}
    \mathcal{P}_{\mathrm{rot}}^\omega(x,t)
      &=(x_1\cos{(t\omega)}-x_2\sin{(t\omega)},x_{2}\cos{(t\omega)}+x_1\sin{(t\omega)},x_3).\label{eq:DefRotationalMovement}
\end{align} to be the rotational movement with angular velocity $\omega>0$. In particular, for any rotation function $\mathcal{P}(x,t)$, with particle $(X_k,V_k,T_k)_k$ define the path $\gamma_k$ as $\mathcal{P}(X_k+(T_k-t)V_k,t)$ for $t\in[T_k,T_{k+1})$.

In the movement defined above, the velocity of the space passing through any fixed point $y\in \mathbb{R}^3$ is constant. For any $\omega>0$, this values creates a constant vector field:
\begin{align}
    \Fr(y)&= (-\omega y_{2},\omega y_1,0)\label{eq:DefRotationalField}.
\end{align} As a physical connection, notice that the derivative in time of the position $D_t [\mathcal{P}_{\mathrm{rot}}^{\omega}(x,t)]$ is equal to $\Fr(\mathcal{P}_{\mathrm{rot}}^{\omega}(x,t))$.

In this paper, we consider three main sources of randomness in the problem. The first source of randomness is the new velocity that the particle has after collision due to the unpredictable angle of contact, the second source is given by the thermal energy in the space that is also transferred to the particles in collision, and the last source of randomness is the time at which the particle interacts with the medium.  Such random properties can depend on the position, in time, and also may have intricate dependencies between them. Here, we are going to ignore such dependencies and define the following independent auxiliary random variables:
\begin{enumerate}
    \item Let $(\eta_k)_{k\in \N}$ be an i.i.d. sequence with a non empty interior support over the open ball $B(0,1)$.
     \item Let $(\Delta_k)_{k\in \N}$ be an i.i.d. sequence of random variables in $\R^3$ distributed as $\mathcal{N}(0,\sigma^2)$ for some small $\sigma>0$. 
    \item Let $(\xi_k)_{k\in \N}$ be an i.i.d. sequence of exponential random variables with rate $\lambda>0$, this is $\P{\xi_k>s}=e^{-s\lambda}$. 
\end{enumerate}
Some considerations over the physical world are made in this choice, see Remark \ref{Rmk:1} for more details.

The variable $\eta$ represents the probability distribution of the final velocity of a particle in a given collision. This collision is specifically a collision between a particle $A$ with an initial velocity of $(1, 0, 0)$ and a mass $m_1$ into a stationary particle $B$ of mass $m_2$ at position $(0, 0, 0)$, where $\eta$ corresponds to the final velocity of the particle $A$. In the general case, with a incident particle with velocity ${v}$ colliding with a target with a speed ${u}$, by changing the referential and normalize the vectors, we define the result of the collision to be:
\begin{align}\label{eq:Rotaçãonacolisão}
  \|{v}-{u}\|R_{{v}-{u}}\cdot \eta +{u} 
\end{align}Where for any vector ${u}$, $R_{{u}}$ is defined to be the three by three rotational matrix that sends $(1,0,0)$ to $u/\|u\|$, preserving angles and orientation.

The variables $(\xi_k)_{k\in \N}$ correspond to the total distance that a particle walks between collisions. This value, together with some calculations, indirectly informs the new time in which the particle interacts with space. To illustrate, consider $(X_k,V_k,T_k)$ a particle in the $k-th$ collision, define the distance walk by the particle till time $t>0$ to be the function:
\begin{align}\label{eq:defdistancerot}
    D_{X_k,V_k}^{\omega}(t)=\int_0^t \left\|V_k-\Fr \left(X_k+sV_k\right)\right\|_{} ds,
\end{align} and with that define $T_{k+1}=T_k+[D^{\omega}_{X_k,V_k}]^{-1}(\xi_k)$, also notice that by construction $|\gamma_k|=\xi_k$.
 
The thermal energy of the space is considered to be an independent velocity $\Delta$ associated with the components of the medium. Such value in this text is used in the collision equation \eqref{eq:Rotaçãonacolisão},  implying that the final velocity of the collision in each collision have the form:
\begin{align}\label{eq:Colisionequation}
    \|{v}-{(u+\Delta)}\|R_{{v}-{(u+\Delta)}}\cdot \eta +{(u +\Delta)}.
\end{align}

\begin{Rmk}\label{Rmk:1}
    The choice of these random variables is based on some assumptions in the real world. It was done in this setting to simplify some of the calculations and provide a reasonable explanation for what should be the phenomenon of collisions.
    
    The collision of particles follows some general Newtonian principles in order to conserve energy and momentum. This leads to two main points: The final velocity after the collision is directly proportional to the relative velocity between the particles; the second is the rotation invariance of the collisions. However, despite the simple calculations, the principal importance of the definition of $\eta$ is the support in the open ball $B(0,1)$, that is, the incident particle always loses energy. In general, we can substitute $\eta$ for any angle distribution that loses a fraction of its initial energy in collisions. 

    With the thermal energy at each point in the space, the relative velocity of a particle becomes a complex definition. We define the relative velocity of a particle as the limit of the average difference between the velocity of the particle and the velocity of a small neighbor around the moving particle. Since $\Delta \sim \mathcal{N}(0,\sigma^2)$, the average relative to $\Delta$ is zero, and this limit is equal to the classical relative velocity. We avoid integrating the norm of the true relative velocity since the integral of $\| \Delta \|$ is not well defined. 

    Since the particle walks into a homogeneous medium, our model supposes that equal lengths in the path imply equal probability of collision. This, together with the loss of memory of the problem, implies that $\xi\sim \mathrm{exp}(\lambda)$ is a reasonable choice of distance. 

\end{Rmk}

Finally, we can state out the main theorems that expose the behaviors of the random variable $X_t$. The first theorem informs us that the particles diverge and follow the velocity of the space:

\begin{Thm}\label{Thm:1}
For every $\omega>0$ the particle $(X_t,V_t,t)$ in a rotational environment satisfies:
\begin{align*}
    d_r(X_t)\to \infty&, \, \text{a.s. and in probability}.\\
    \frac{V_t}{\Fr(X_t)}\to 1&, \, \text{a.s. and in probability}.
\end{align*}Moreover, there exists an value $\beta>0$ such that:
\begin{align*}
    \P{\liminf\limits_{n\to \infty}{\{d_r(X_n)>\beta n\}}}=1,
\end{align*} where for sets $(A_n)_n$, $\liminf\limits_{n\to \infty} A_n=\bigcup\limits_{n=1}^{\infty} \bigcap\limits_{m>n} A_m$.  
\end{Thm}

In particular, for some $t_0>0$, for every $t>t_0$ the radial position of the particle $X_t$ grows better than some linear function over time. Also, since the space is rotating and $V_t/F(X_t)$ converges to one, the particle is performing a spiral movement.

The spiral movement can be described using a decomposition of velocities in tree coordinates that depends on the position of the particle. The decomposition will be in the height, radial, and azimuth components, $(v^h,v^r,v^a)$, respectively. Given a particle in $x$ with velocity $v$, we define the height component $v^h$ as the velocity of the third coordinate position, the radial component $v^r$ will be the component in the subspace generated by $\langle(x_1,x_2,0)\rangle=\{s(x_1,x_2,0): s\in \R\}$, and the last component $v^a$ will belong to the tangent of the subspace generated by $\langle v^r,v^h\rangle=\{av^r+bv^h, a,b\in \R\}$. 

Due to collisions, the particle does not take a continuous path and neither does it have a continuous transition of the velocity. To define a value that we can associate with velocity, one should observe a conditional expectation of the discrete version of velocity given the distance of the particle. 
Then, at the $k-$th collision consider the difference of the radial distances as $d_r(X_{k+1})-d_r(X_k)$,  and define the radial velocity as the random variable: 
\begin{align*}
    V^r_k=\Ex{\frac{d_r(X_{k+1})-d_r(X_k)}{T_{k+1}-T_k}\middle| d_r(X_k)}
\end{align*}That represent the best random value of the ratio $\frac{d_r(X_{k+1})-d_r(X_k)}{T_{k+1}-T_k}$, given the position of the particle $d_r(X_k)$. 

The second theorem will bound above the movement $X_t$ proving a sharp bound on the expected radial velocity. Showing that Langevin's approximation in equation \eqref{eq:Intro-Lagevin-Possible-Approximation} is not true; the particle indeed receives a force to diverge, but this force does not grow as the centripetal force. 
 
\begin{Thm}\label{Thm:2}
    There exists a radial distance $N_0>0$, such that when $d_r(X_k)>N_0$ the time between collisions is shorter than $\frac{2\max\{\xi,2\}}{\omega\sqrt{d_r(X_k)}}$ and bigger than $\sqrt{\frac{\xi}{\omega^2d_r(X_k)}}$. Furthermore, there exist constants $\infty>c_2>c_1>0$ such that: 
    \begin{align*}
        \P{\liminf\limits_{n\to \infty} \left\{\frac{V^r_n}{\sqrt{n}}\in(c_1,c_2)\right\}}=1. 
    \end{align*}
\end{Thm}

\begin{Rmk}\label{Rmk:Phiscs}
   The particle at $(r,0,0)$ with the same velocity of space takes a constant amount of time to go around the circle with radius $r$. In this constant time, the particle will drift away a distance of order $\sqrt{r}$, very small comparable to $r$. As the path of the particle resembles the circle, a centripetal force should appear. However, this evaluation is wrong: In this problem, despite the collisions, the particle walks almost everywhere in a straight line, obeying every physic law. In collisions, the conservation of momentum and energy is due to the loss of energy. Thus our problem is obeying every Newtonian law. 
   
   The result of Theorem \ref{Thm:2} should not be understood from a continuous physics perspective, but rather in the random sense, where several collisions result in a strong law that determines its movement. 

   The physical stipulation done in Section \ref{sec:Intro} takes two limits at the same time. The first limit is related to the rotation velocity, and the second limit is related to the centripetal force of order $r$, both values can be seen in the Lavengin's Equation \eqref{eq:Intro-Lagevin-Possible-Approximation}. As Theorem \ref{Thm:2}, the centripetal force of the problem has order $\partial_r\sqrt{r}=\frac{1}{2\sqrt{r}}$, showing that we cannot assume that the limits can be taken in the problem, and for the authors making the result interesting. 
\end{Rmk}

%% file: 3.tex
\section{The Rotational Environment}\label{sec:rot}
\noindent

This Section is devoted to expose the behavior of a random walk in a medium that rotates with some angular velocity $\omega>0$, such movement is defined in equation \eqref{eq:DefRotationalMovement} in Section \ref{sec:Notation}. In particular, we fixed $\omega>0$, $\lambda>0$, and $\sigma>0$. Let $(X_k, V_k, T_k)_{k\in \N}$ be the position and velocity of the particle over the colliding times $(T_k)_k$, and also fix the distributions of the sequences $(\xi_k)_k$, $(\eta_k)_k$, and $(\Delta_k)_k$. 

To control the position $X_t$, the first step is to understand the behavior of the vector field $\Fr$ between collisions. This value indirectly will be fundamental to bound the fluctuations of the expected velocity in Proposition \ref{Prop:1}. 

\begin{Lem}\label{Lem:1}
    Fixed the sequence $(\xi_k)_k$ of distances walked between collisions, and let $(X_k)_k$ be the set of positions in the interaction times $(T_k)_k$. For every $\omega > 0$, we have:
    \begin{align*}
        |\Fr(X_k) - \Fr(X_{k+1})| < \omega \xi_k.
    \end{align*}
\end{Lem}
\begin{proof}[Proof of Lemma \ref{Lem:1}]
    Fix a constant $\xi > 0$, and consider for any $z \in \mathbb{R}$ the point $p_z = (z, 0, 0)$. Then, define the quantity: 
    \begin{align*}
        N = \sup\limits_{z \in \mathbb{R}} \max\limits_{x \in B(p_z, \xi)} |\Fr(p_z) - \Fr(x)|.
    \end{align*}Using the fact that $|\Fr(p_z) - \Fr(x)|$ is a continuous function in $x$ with constant derivatives towards the axis, the Lagrangian method shows that the maximum difference of this is achieved at the points $p_z \pm (\xi, 0, 0)$. Substituting these values, we get $N = \omega \xi$.

     Since the space is Euclidean and the movement is continuous, for any $t_k < t < t_{k+1}$, we have $X_t \in B(\mathcal{P}^\omega_{\text{rot}}(X_k, t), \xi_k)$, and thus:
    \begin{align*}
        |\Fr(X_k) - \Fr(X_{k+1})| = |\Fr(\mathcal{P}^\omega_{\text{rot}}(X_k, t_{k+1})) - \Fr(X_{k+1})| \leq N = \omega \xi_k,
    \end{align*}as desired.
    
\end{proof}

The particle is expected to align with the movement of the medium. To control such alignment define the difference between the velocity in the $k$-th collision and the velocity of the field as follows:
\begin{align}\label{eq:fluctuations}
    E_k = V_k - \Fr(X_k),
\end{align}such value is interpreted as an error or a fluctuation. In particular, writing $V_k = \Fr(X_k) + E_k$ the velocity of the particle becomes the velocity of the medium plus some fluctuation. The goal in the next proposition is to show that the random variables $(E_k)_k$ have a uniformly bounded mean, implying that when $\Fr(X_k)$ diverges, the velocity of the particle becomes the velocity of the space with high probability.

\begin{Prop}\label{Prop:1}
    For simplicity, in this proposition, let $\eta = E(\| \eta \|) \in (0, 1)$ and $\Delta = E(\| \Delta \|) > 0$. Then, we have for a particle with $X_0=x_0$ and $V_0=v_0$:
    \begin{align*}
        \Ex{\|E_k\|} <  \frac{\Delta \lambda+\omega \eta}{\lambda(1 - \eta)}+ \|v_0-\Fr(x_0)\|\eta^k.
    \end{align*}
\end{Prop}

\begin{proof}[Proof of Proposition \ref{Prop:1}] The proof is going to be  an induction of the mean value of the error, where we can associate it by the last collision, take a note on the importance of the loose of energy in the collisions, i.e. $\eta<1$. By computing the collision equation \eqref{eq:Colisionequation} for the particle $(X_k,V_k,t_k)$, we get : 
    \begin{align*}
        V_k = \|V_{k-1} - (\Fr(X_k) + \Delta_k)\| R_{V_{k-1} - (\Fr(X_k) + \Delta_k )} \cdot \eta_k + (\Fr(X_k) + \Delta_k).
    \end{align*}
    Thus, by the triangular inequality, and the fact that the rotational matrix $R$ has norm one, one can get that:
    \begin{align*}
    \|V_k& - \Fr(X_k)\| \leq \|V_{k-1} - (\Fr(X_k) + \Delta_k)\|\|\eta_k\| + \|\Delta_k\|\\
    \|E_k\|&\leq (\|E_{k-1}\|+\|\Fr(X_{k-1})-\Fr(X_k)\|+\|\Delta_k\|)\|\eta_k\|+\|\Delta_k\|.
    \end{align*} Applying Lemma \ref{Lem:1}, we get the induction equation that relates $\|E_k\|$ with $\|E_{k-1}\|$:
    \begin{align}\label{eq:inequalityovertheerror}
        \|E_k\| \leq (\|E_{k-1}\| + \omega \xi_k + \|\Delta_k\|) \|\eta_k\| + \|\Delta_k\|.
    \end{align}Therefore, by looking into the expectation, and by using independence between the random variables, we get by induction:
    \begin{align*}
       \Ex{\|E_k\|} &\leq  (\Ex{\|E_{k-1}\|} + \frac{\omega}{\lambda} + \Delta) \eta + \Delta\\
       &\leq ((\Ex{\|E_{k-2}\|} + \frac{\omega}{\lambda} + \Delta) \eta + \Delta + \frac{\omega}{\lambda} + \Delta) \eta + \Delta,\\
       &\leq \Ex{\|E_0\|}\eta^k+\sum_{\ell=1}^k \left(\Delta + \frac{\omega}{\lambda} \eta\right) \eta^{\ell-1}\\
       &<\|v_0-\Fr(x_0)\|\eta^k+ \frac{\Delta \lambda+\omega \eta}{\lambda(1 - \eta)},
    \end{align*} as desired.
\end{proof}

\begin{Rmk} \label{Rmk:Prop1}
By just applying the square to both sides of the equation \eqref{eq:inequalityovertheerror}, it is also possible to bound the second moment of $\|E_k\|$. Using the second moment, one can make the proof a little sharper, but it is an effort that does not creates a greater contribution to the paper.  

Also, observe that if there exists $\e>0$ such that $\|\eta_k\|<1-\e$, a similar inequality appears by having an exponential decay with rate $\ln{(1-\e)}$. The unique point used to bound the fluctuation is the loss of energy in the collisions.

\end{Rmk}\bigskip


By Proposition \ref{Prop:1}, $\Ex{\|E_k\|}$ is bounded for every $k$, and exponentially fast can be bounded by a reasonable number that does not depend on the initial condition. This indicates that the fluctuation between the velocity of the particle and the velocity of the space is with high probability close to the value of $\Ex{\|E_k\|}$. This fact becomes stronger when $d_r(X_k)$ increases, forcing the error value to become smaller compared to the value of $\Fr(X_k)$. In particular, the velocity of the particle aligns with the velocity of the medium whenever the particle diverges. To control this alignment in a sharper result, observe the following geometrical lemma: 

\begin{Lem}\label{Lem:2}
    Consider a particle at position $x$ with velocity $v$, such that $v = \Fr(x) + e$, where $ e = (e_1, e_2, e_3)$. Let $C>0$ be any constant, then there exists $N_1 = N_1 (C, \omega, \|e\|)$ such that for all $x$ where $d_r(x) > N_1$ the total distance traveled by the particle and its position radius satisfy at time $\frac{C}{\sqrt{d_r(x)}}$:
    \begin{align*}
        D_{x,v}^{\omega}\left(\frac{C}{\sqrt{d_r(x)}}\right)\leq \omega^2C^2,\\
        d_r\left(x+\frac{C}{\sqrt{d_r(x)}}v\right)- d_r(x)\geq   \frac{\omega^2 C^2 }{6}.
    \end{align*} Moreover, for every $t>\frac{C}{\sqrt{d_r(t)}}$, we also gets:
    \begin{align}\label{eq:tgrandeoraiocresce}
         d_r\left(x+tv\right)-d_r(x)\geq   \frac{\omega^2 C^2 }{6}.
    \end{align} Forcing a increase of radius for the particle.
\end{Lem}

\begin{proof}[Proof of Lemma \ref{Lem:2}]
Fixing the values of $\|e\|>0$, $C>0$ and $\omega>0$, one can find $N_1=N_1(C,\omega,\|e\|)$ such that, for all $n>N_1$:
\begin{align*}
    &\frac{\omega^2 C^2}{2n} < 1,\\
        &\frac{C\|e\|}{\sqrt{n}}\left(1+ \frac{C\omega}{2\sqrt{n}}\right)<\frac{\omega^2  C^2}{2},\\
        &\left| \frac{2Ce_1}{\sqrt{n}}+\frac{2C^2\omega e_2}{n}+\frac{C^2d_r^2(e)}{n^2}\right|<\frac{\omega^2 C^2}{2}.
\end{align*}

Observe in Figure \ref{fig:1} all the quantities related to the lemma. The challenge of the lemma is associated with the case when the particle wants to go in the direction of the center; in this case, since the velocity of the particle is aligned with the velocity of the medium, the path that the particle does gets out of the circle in a short period of time due to its high velocity. Moreover, the distance walked in this path is also small. 
\begin{figure}
    \centering
    \includegraphics[scale=1]{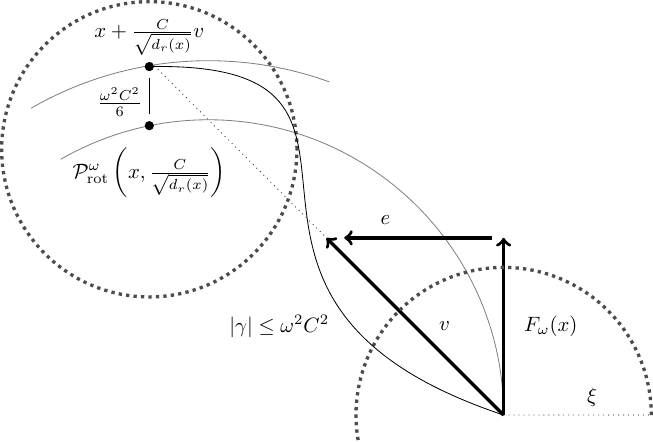}
    \caption{A particle in the position $x$ walking in direction $v=\Fr(x)+e$. At time $\frac{C}{\sqrt{d_r(x)}}$ the particle already get out the circle with norm $d_r(x)$ and it is located around the circle with radius $d_r(x)+ \frac{\omega^2 C^2}{6}$. The path in green have size at most $\omega^2 C^2$, and it never gets out of a distance of $\xi$ of the circle. }
    \label{fig:1}
\end{figure}

To start proving, fix $x$ such that $d_r(x)>N_1$ and notice that $\Fr$ is a linear function. Using the triangular inequality and the fact that $\|\Fr(x)\|=\omega d_r(x)\leq \omega \|x\|$,  one gets:  
    \begin{align*}
          D_{x,v}^{\omega}(t)&=\int_0^t \|v - \Fr(x + v s)\| ds = \int_0^t \|v - \Fr(x) + s\Fr(v)\| ds \\
          &=\int_0^t \|e+ s\Fr(\Fr(x))+ s\Fr(e)\| ds\leq \|e\| t+ \frac{\omega t^2}{2} d_r(e)+ \frac{\omega t^2}{2} d_r(x)\\&\leq \|e\|\left(t+\frac{\omega t^2}{2}\right) + \frac{\omega^2 t^2}{2}d_r(x).
    \end{align*} Thus, since $d_r(x)>N_1$, we get:
    \begin{align*}
         D_{x,v}^{\omega}\left(\frac{C}{\sqrt{d_r(x)}}\right)&\leq \frac{C\|e\|}{\sqrt{d_r(x)}}\left(1 + \frac{\omega C}{2\sqrt{d_r(x)}} \right) + \frac{\omega^2 C^2}{2}\leq \omega^2 C^2.
    \end{align*}As desired.
    
    To understand the position of the particle at time $\frac{C}{\sqrt{d_r(x)}}$, we can perform a direct calculation. However, using symmetries of the problem by rotating the space or changing variables, the computation becomes shorter. So, assume without loss of generality that $x = (z, 0, 0)$ for some $z > 0$, in particular $d_r(x) = z$, then:
\begin{align*}
    d_r\left(x + \frac{C}{\sqrt{d_r(x)}}v\right)&= d_r\left(z+\frac{Ce_1}{\sqrt{z}},\frac{Ce_2+C\omega z}{\sqrt{z}},\frac{Ce_3}{\sqrt{z}}\right)\\
    &=\sqrt{\left(z+\frac{Ce_1}{\sqrt{z}}\right)^2+ \left(\frac{Ce_2+C\omega z}{\sqrt{z}}\right)^2}\\
    &=\sqrt{ z^2+ 2Ce_1 \sqrt{z}+ \frac{(Ce_1)^2}{z} + \frac{(Ce_2)^2}{z} + 2C^2\omega e_2+ (\omega C)^2 z}\\
    &\geq \sqrt{z^2+z\left(\omega^2C^2-\left|\frac{2Ce_1}{\sqrt{z}}+\frac{2C^2\omega e_2}{z}+\frac{C^2d_r^2(e)}{z^2}\right|\right)}\\
    &\geq \sqrt{z^2 + z\frac{\omega^2 C^2}{2}}
\end{align*} Using that $d_r(x)>N_1$ one gets that $\frac{C^2\omega^2}{2d_r(x)}<1$. Thus, since $\sqrt{1 + x} > 1 + \frac{x}{3}$ for every $x\in(0,1)$, we get for $d_r(x)>N_1$ that:
\begin{align*}
     d_r\left(x + \frac{C}{\sqrt{d_r(x)}}v\right)\geq  d_r(x) +\frac{\omega^2 C^2}{6}.
\end{align*}

Notice that a particle starting with radius $d_r(x)$ and moving in a line in one direction can hit a cylinder with radius larger than $d_r(x)$ at most at one point. And whenever it hits this point, as time passes, the distance between the point and the cylinder increases. This proves equation \eqref{eq:tgrandeoraiocresce}, finishing the proof.
\end{proof}

The result of Lemma \ref{Lem:2} implies that by simply using a bound in the fluctuation, suggested in Proposition \ref{Prop:1}, the movement of the particle is somehow deterministic, since it cannot go back to near the origin with high probability. This blockage created by a combination of geometry and randomness generates a virtual force on the particle, it is virtual since it does not fundamentally exist and it is the result of many collisions in this faraway point. 

The value, module, and behavior of this virtual force created by a sum of factors can be studied by comparing the movement of the particle with a bias random walk. There exist many instruments that can control the bias random walk , see \cite{Durrett}, here we state a version of Cramer's Theorem:

\begin{Thm}[Crammer]\label{Thm:Crammer}
    Let $(Y_k)_k$ be a sequence of independent random variables, each distributed as $Y$, and let $S_n=Y_1+\cdots + Y_n$. Let $\Lambda(\theta)= \log{\Ex{e^{\theta Y}}}$, and suppose that there exists $\delta>0$ such that for every $\theta\in (-\delta,\infty)$ we have $\Lambda(\theta)<\infty$. Define $\Lambda^*(x)=\sup\limits_{\theta\in (-\delta,\infty)} \{\theta x-\Lambda(\theta)\}$, and for every $\e>0$ define the function $I(\e)=\inf\left\{\Lambda^*(x):|x-\Ex{Y}|>\e\right\}$. Assuming that $I(\e)> 0$ for every $\e>0$,we get
    \begin{align*}
        \P{\left|\frac{S_n}{n}-\Ex{Y}\right|> \e}\leq e^{-nI(\e)}. 
    \end{align*}An thus, we have
    \begin{align}\label{eq:forallnSnbehaved}
        \P{\bigcup_{n}\left\{ \left|\frac{S_n}{n}-\Ex{Y}\right|> \e\right\}}>0. 
    \end{align} 
\end{Thm}

The proof of this Theorem is an adaptation of the famous Crammer Theorem presented in \cite{Durrett}. The choice to use this theorem is based on the last affirmation, where by the exponential bound on the events, one can use Borel-Cantelli Theorem to get the equation~\eqref{eq:forallnSnbehaved}.

With this result, one can prove Theorem \ref{Thm:1}, and get an overall behavior of the random variable $X_t$.

\begin{proof}[Proof of Theorem \ref{Thm:1}] 

The proof is based on the coupling between a bias random walk and the movement of the particle. By Lemma \ref{Lem:2}, particles with a certain distance from the origin have in each collision a high probability chance of being dragged further away. This affirmation that becomes stronger as the distance of the particle grows makes the transition between the randomness of the problem into a deterministic movement. 

The position of the particle is random, but its limit behavior has tree possibilities: The first possibility is when the particle does not leave a certain bounded region, the second case is when the particle diverges and never returns to near the origin, and the last possibility is when it does excursions going to faraway points and returning.  

The particle could not be trapped in a finite region. For each bounded region, notice that a particle with initial velocity $v$ after a collision with other at velocity $u$, have a final velocity $v_f$ equal to:  
\begin{align*}
    v_f=\|{v}-{(u+\Delta)}\|R_{{v}-{(u+\Delta)}}\cdot \eta +{(u +\Delta)}.
\end{align*}Since $\Delta$ is a continuous random variable, we get that $v_f$ is zero with probability zero, and with bounded values of $u$, and $\eta$, we get $\Ex{\|v_f\|}\geq \Ex{\|\Delta\|}>0$. Thus with some positive probability, the velocity in a bounded region achieve some non negligible values. With any value fixed of velocity, one can find a distance $D$ such that for every point inside the bounded region, walking a distance $D$, we get out, since $\{\xi>D\}$ is possible with positive probability. By Borel Cantelli argument, the particle gets out the circle. 

To finish the argument, it rests to show that the particles diverge and never return. For this, let $\Delta=\Ex{\|\Delta\|}$, $\eta=\Ex{\|\eta\|}$, fix $\lambda>0$, and $\omega>0$. Define the constants $C$ and $p\in(0,1)$ such that:
\begin{align*}
    C=\sqrt{\frac{1}{7\lambda\omega^2}},\text{ and }
    p=\frac{1}{14} (7 - 6 e^{1/7}).
\end{align*} 

By Lemma \ref{Lem:2}, find $k_0=k_0(x_0,v_0,\eta)$, such that for every $k>k_0$, we get that:
\begin{align*}
        \Ex{\|E_k\|}\leq 2\frac{\Delta\lambda+\omega \eta}{\lambda(1-\eta)}. 
\end{align*} Thus, with the fixed value of $p$, by the Markov inequality one may get for every $k>k_0$, that:
\begin{align*}
        \P{\|E_k\|\geq 4\frac{\Delta\lambda+\omega \eta}{p\lambda(1-\eta)}}\leq p,\text{ and }
        \P{\xi_k\geq \omega^2 C^2}=e^{-\lambda C^2 \omega^2}.
\end{align*} 

In each collision of the particle, define the events: 
\begin{align*}
    G_k&=\left\{\xi_{k+1}\geq\omega^2C^2 , \|E_k\|< 4\frac{\Delta\lambda+\omega \eta}{p\lambda(1-\eta) }\right\}\\
    B_k^1&=  \left\{\xi_{k+1}<\omega^2C^2\right\} \\
    B_k^2&= \left\{\xi_{k+1}>\omega^2C^2 , \|E_k\|\geq 4\frac{\Delta\lambda+\omega \eta}{p\lambda(1-\eta) }\right\}
\end{align*} And, with that define the random variable:
\begin{align*}
    Y_k= \frac{\omega^2 C^2}{6} \ind\{G_k\} -\xi_{k+1}\left(\ind\{B_k^1\}+\ind\{B_k^2\}\right)
\end{align*}

Using Lemma \ref{Lem:2}, find $N_1=N_1\left(C,\omega,2\frac{\Delta\lambda+\omega \eta}{p\lambda(1-\eta)}\right)$, such that whenever $d_r(X_k)>N_1$ and $G_k$ occurs, the total distance walked by the particle and its position radius satisfy at time $T_k+\frac{C}{\sqrt{d_r(X_k)}}$, that:
\begin{align*}
        D_{X_k,V_k}^{\omega}\left(\frac{C}{\sqrt{d_r(X_k)}}\right)\leq \omega^2C^2,\text{ and }
        d_r\left(X_k+\frac{C}{\sqrt{d_r(X_k)}}V_k\right)- d_r(X_k)\geq   \frac{\omega^2 C^2 }{6}.
\end{align*} Moreover whenever $G_k$ does not happens, one may assume that $d_r(X_{k+1})>d_{r}(X_k)-\xi_k$. Therefore, whenever $d_r(X_k)>N_1$, we get that $Y_k<d_r(X_{k+1})-d_r(X_{k})$, implying that the summation of $Y_k$ starting in $N_1$ when $d_r(X_k)>N_1$ bounded bellow the radii value of $X_t$. If $\sum_k Y_k$ eventually becomes zero, we cannot guarantee that the particle is indeed greater than $N_1$, but in this case, one can use the same Borel Cantelli argument that shows that a particle never stays in 
 a bounded region to finally reach $N_1$ again and repeat the summation of $Y_k$. Now, if eventually the $\sum_k Y_k$ diverges, we have that $(d_r(X_k))_k$ diverges together. To show that the summation can diverge with positive probability, let us use Cramer's theorem \ref{Thm:Crammer}. For this, let $B\sim \mathrm{Ber}(p)$, and define:
\begin{align*}
    Y= \frac{\omega^2 C^2}{6} \ind\{\xi_k>\omega^2 C^2, B=0 \} -\xi_{k+1}\left(\ind\{\xi_k\leq\omega^2 C^2\}+\ind\{\xi_k>\omega^2 C^2,B=1\}\right).
\end{align*} Notice that by a simple coupling of the random variables $B$, one should get that $Y\leq Y_k$. Therefore, comparing $\sum_k Y_k$ with a sum of independent copies of the random variable $Y$, we are in the context of Crammers theorem, and:
\begin{align*}
    \Ex{Y} &= \frac{\omega^2 C^2}{6}(1-p)e^{-\lambda \omega^2C^2}-\int_0^{\omega^2 C^2} x \lambda e^{-\lambda x} dx- p\int_{\omega^2 C^2}^{\infty} x\lambda e^{-\lambda x} dx\\
    &= \frac{\omega^2 C^2}{6}(1-p)e^{-\lambda \omega^2C^2}-\frac{1 - e^{-\lambda \omega^2 C^2}(1 +  \lambda \omega^2 C^2)}{\lambda} -p \frac{ e^{-\lambda \omega^2 C^2}(1 +  \lambda \omega^2 C^2)}{\lambda}\\
    &=\frac{1}{6\lambda} \left(-6 + e^{-\lambda \omega^2 C^2} (1 - p) (6 + 7\lambda C^2 \omega^2)\right).
\end{align*} By the choice of $C$ and $p$, we get:
\begin{align*}
    \Ex{Y} &\geq \frac{1}{6\lambda} \left(-6 + 7 e^{-1/7} (1 - p)\right)>0 .
\end{align*}

To finish the proof, let us show that the random variable $Y$ has a positive rate function $I(\e)$ for every $\e>0$. For this, let $\theta>-\lambda$, thus: 
\begin{align}
    \Ex{e^{\theta Y}}
    &=(1-p)\left(e^{\omega^2 C^2 \left(\frac{\theta}{6}-\lambda\right)}- \frac{\lambda e^{-\omega^2 C^2(\theta+\lambda)}}{\theta+\lambda}\right)+\frac{\lambda}{\theta+\lambda}. \label{eq:ExpthetaYfunct}
\end{align} In particular, for every fixed value of $p,C,\lambda,$ and $\omega$, the function in equation \eqref{eq:ExpthetaYfunct} is analytic in $\theta$, and not equal to some purely exponential function. Therefore, $\Lambda^*(x)>0$ uniformly  for every $x$ not close to $\Ex{Y}$. Then, for every  $\e>0$,  we get $I(\e)>0$, and  we have that by Theorem \ref{Thm:Crammer}: 
\begin{align*}
    \P{\left|\frac{S_n}{n}-\Ex{Y}\right|>\e,\forall n}>0.
\end{align*} Where $S_n=\sum_{k=1}^n Y^{(k)}$ for i.i.d. random variables $Y^{(k)}$ distributed as $Y$. Thus with positive probability the particle when reach an radii distance of $N_1$ diverges.

Since by Borel Cantelli it will always return, we get with probability one the particle diverges, thus proving that almost surely and in probability we have that:
\begin{align*}
    d_r(X_t)\to \infty.
\end{align*} Then, we get by Proposition \ref{Prop:1}, that almost surely and in probability we get:
\begin{align*}
    \frac{V_t}{\Fr(X_t)}\to 1. 
\end{align*} 

More than this, the exponential decay of Crammer's Theorem together with the Borel-Cantelli argument shows that for every $\e>0$, eventually we get $d_r(X_n)>n(\Ex{Y}-\e)$, for every $n$. By choosing $\e>0$ accordingly, we finish the proof as desired. 
\end{proof}

To prove Theorem \ref{Thm:2}, a procedure similar to Lemma \ref{Lem:2} will be done. Now, we are going to give an upper bound on the distances walked by the particle, in relation to the time. Therefore, consider: 
\begin{Lem}\label{Lem:3}
    Consider a particle at position $x$ with velocity $v$, such that $v = \Fr(x) + e$, where $ e = (e_1, e_2, e_3)$.  There exists a distance $N = N (\omega, \|e\|)$ such that for all $x$ where $d_r(x) > N$, the distance traveled by the particle in time $\frac{2\max\{2,\xi\}}{\omega\sqrt{d_r(x)}}$ satisfies:
    \begin{align*}
        D_{x,v}^{\omega}\left(\frac{2\max\{2,\xi\}}{\omega\sqrt{d_r(x)}}\right)\geq \xi.
    \end{align*} 
\end{Lem}
\begin{proof}[Proof of Lemma \ref{Lem:3}]
Fixing the values of $\|e\|>0$, and $\omega>0$, one can find $N_2=N_2(\omega,\|e\|)$ such that, for all $n>N_2$:
\begin{align*}
    &\frac{2d_r(e)}{n} < \frac{\omega^2}{2}\text{ and }\frac{4\|e\|^2\omega^2}{n} <1.
\end{align*}

To start proving, fix $x$ so that $d_r(x)>N_2$ and notice that $\Fr$ is a linear function. Using the triangular inequality and the fact that $\|\Fr(x)\|=\omega d_r(x)\leq \omega \|x\|$, one gets
\begin{align*}
          D_{x,v}^{\omega}\left(\frac{2\max\{2,\xi\}}{\omega\sqrt{d_r(x)}}\right)&=\int\limits_0^{\frac{2\max\{2,\xi\}}{\omega\sqrt{d_r(x)}}} \|e+ s\Fr(\Fr(x))+ s\Fr(e)\| ds\\
          &\geq \int\limits_{\frac{4\|e\|}{d_r(x)}}^{\frac{2\max\{2,\xi\}}{\omega\sqrt{d_r(x)}}} \|e+ s\Fr(\Fr(x))+ s\Fr(e)\| ds
    \end{align*} Since, for every $s> \frac{4\|e\|}{d_r(x)}$, we get $s\omega^2 d_r(x)\geq  2s d_r(e)+2\|e\|\geq 2\|s\Fr(e)+e\|$, then by the reverse triangular inequality, assuming that $d_r(x)>N_2$ one gets that:
    \begin{align*}
          D_{x,v}^{\omega}\left(\frac{2\max\{2,\xi\}}{\omega\sqrt{d_r(x)}}\right)
          &\geq \int\limits_{\frac{4\|e\|}{d_r(x)}}^{\frac{2\max\{2,\xi\}}{\omega\sqrt{d_r(x)}}} \frac{1}{2}\| s\Fr(\Fr(x))\| ds= \left(\max\{2,\xi\}\right)^2 - \frac{4\|e\|^2\omega^2}{d_r(x)} 
          \geq \xi.
    \end{align*} Where if $\xi>2$, $\xi^2-1>\xi$, and if $\xi<2$, it is also clear that $\xi<3$.  
\end{proof}

Finally, we are able to proof Theorem \ref{Thm:2}, and find the radial velocity of the particle. 

\begin{proof}[Proof of Theorem \ref{Thm:2}] The proof of Theorem \ref{Thm:1} gives an non trivial lower bound for the positions of $d_r(X_k)$ when $k$ is big. To give an upper bound trivially, one may use Lemma \ref{Lem:1} to get that $|d_r(X_{k+1})-d_r(X_{k})|<\omega \xi_{k}$. In that way, using the same Crammer's theorem now in each step walking faraway from the origin a distance $\xi_k$, one may get that exists constants $\alpha>\beta>0$ such that:
 \begin{align}\label{eq:wheretheyare}
        \P{\liminf\limits_{n\to \infty} \left\{\beta n< d_r(X_n)<\alpha n\right\}}=1. 
\end{align}

Therefore, given one particle in the rotational environment, there exists a collision $k_0$ such that for all $k>k_0$, we get that $\{\beta k< d_r(X_k)<\alpha k\}$ occurs. To finish the proof, take $C=\sqrt{\frac{\xi}{\omega^2}}$ from Lemma \ref{Lem:2} and fix $N_0=\max\{N_1,N_2\}$, where $N_1$ is the constant of Lemma \ref{Lem:1} and $N_2$ is the constant of Lemma \ref{Lem:2}. Then, whenever $d_r(X_k)>N_0$ and $k>k_0$, the time between collisions is greater than $\frac{C}{\sqrt{d_r(X_t)}}$, and we get that:
\begin{align}
    V_k^r\leq \Ex{\frac{\alpha}{\frac{\sqrt{\xi}}{\omega^2\sqrt{d_r(X_k)}}}\middle | d_r(X_k)}
    =\alpha \omega^2 \sqrt{d_r(X_k)} \Ex{\xi^{-1/2}}
    = \alpha \omega^2 \sqrt{\lambda \pi} \sqrt{d_r(X_k)}\label{eq:LowerboundV}
\end{align} And, for the other side, using Lemma \ref{Lem:3}, the time between collision is smaller than $\frac{2\max\{2,\xi\}}{\omega\sqrt{d_r(X_k)}}$, thus whenever $d_r(X_k)>N_2$ and $k>k_0$, there exists a constant $c>0$ such that:
\begin{align}\label{eq:UpperboundV}
    V_k^r\geq \Ex{\frac{\beta}{\frac{2\max\{2,\xi\}}{\omega\sqrt{d_r(X_k)}}}\middle | d_r(X_k)}
    = \frac{\beta\omega}{2}  \Ex{\min\{1/2, \xi^{-1}\}} \sqrt{d_r(X_k)} =\frac{\beta\omega}{2} c \sqrt{d_r(X_k)}.
\end{align} 

By equation \eqref{eq:wheretheyare}, using the upper bound of equation \eqref{eq:UpperboundV}, and lower bound of equation \eqref{eq:LowerboundV} , one may get that the radial velocity satisfies that:
 \begin{align*}
        \P{\liminf\limits_{n\to \infty} \left\{\frac{V^r_n}{\sqrt{n}}\in\left(\frac{\beta\omega}{2} c,\alpha \omega^2 \sqrt{\lambda \pi} \right)\right\}}=1. 
    \end{align*}As desired. 

\end{proof}

%% file: main.bbl
\begin{thebibliography}{1}

\bibitem{Durrett}
Richard Durrett.
\newblock {\em Probability: theory and examples}.
\newblock Duxbury Press, Belmont, CA, second edition, 1996.

\bibitem{RandonWalkExternalForce1}
Michael J.~Plank Edward A.~Codling and Simon Benhamou.
\newblock Random walk models in biology.
\newblock {\em J. R. Soc. Interface.}, page 813–834, 2008.

\bibitem{RandonWalkExternalForce3}
Gerald~W. {Englert}.
\newblock {Random Walk Theory of Elastic and Inelastic Time Dependent Collisional Processes in an Electric Field}.
\newblock {\em Zeitschrift Naturforschung Teil A}, 26(5):836--848, May 1971.

\bibitem{BrownianExperiment}
Brian~J. Ford.
\newblock Robert brown, brownian movement, and teethmarks on the hatbrim.
\newblock {\em The Microscope. 39: 161–171.-}, 1991.

\bibitem{RandonWalkExternalForce2}
F.~Le~Vot and S.~B. Yuste.
\newblock Continuous-time random walks and fokker-planck equation in expanding media.
\newblock {\em Phys. Rev. E}, 98:042117, Oct 2018.

\bibitem{Rotationparticleslab}
Hartmut L\"owen.
\newblock Active particles in noninertial frames: How to self-propel on a carousel.
\newblock {\em Phys. Rev. E}, 99:062608, Jun 2019.

\bibitem{Patlak}
Clifford~S. Patlak.
\newblock Random walk with persistence and external bias.
\newblock {\em The bulletin of mathematical biophysics}, 15:311--338, 1953.

\bibitem{Langevin-Book}
L.~C.~G. Rogers and David Williams.
\newblock {\em Diffusions, Markov Processes and Martingales}, volume~2 of {\em Cambridge Mathematical Library}.
\newblock Cambridge University Press, 2 edition, 2000.

\bibitem{Langevin-Intr-Apl-Book}
N.G. {Van Kampen}.
\newblock Chapter ix - the langevin approach.
\newblock In N.G. {Van Kampen}, editor, {\em Stochastic Processes in Physics and Chemistry (Third Edition)}, North-Holland Personal Library, pages 219--243. Elsevier, Amsterdam, third edition edition, 2007.

\end{thebibliography}
